\documentclass[12pt,a4paper]{article}

\usepackage[dvips]{graphics}
\usepackage[dvips]{epsfig}

\newtheorem{lemma}{Lemma}

\newtheorem{prop}[lemma]{Proposition}

\newcommand{\dimo}{\noindent \emph{Proof. }}
\newcommand{\qed}{\\ \rightline{$\Box$}\\}

\usepackage{latexsym}
\usepackage{amsfonts}
\usepackage{amssymb}
\usepackage{amsmath}
\usepackage[english]{babel}
\begin{document}

\title{COMPUTING MATVEEV'S COMPLEXITY \\ VIA CRYSTALLIZATION
THEORY:\\
THE ORIENTABLE CASE  \footnote{Work performed under the auspicies
of the G.N.S.A.G.A. of the C.N.R. (National Research Council of
Italy) and financially supported by M.U.R.S.T. of Italy (project
``Strutture geometriche delle variet\`a reali e complesse'') and
by Universit\`a degli Studi di Modena e Reggio Emilia (project
``Strutture finite e modelli discreti di strutture geometriche
continue'').}}
\author{Maria Rita CASALI  - Paola CRISTOFORI \\
Dipartimento di Matematica Pura ed Applicata \\
       Universit\`a di Modena e Reggio Emilia \\
       Via Campi 213 B \\  I-41100 MODENA (Italy)}

\maketitle

\begin{abstract}
\par \noindent
By means of a slight modification of the notion of {\it
GM-complexity} introduced in \cite{[C$_4$]}, the present paper
performs a graph-theoretical approach to the computation of {\it
(Matveev's) complexity} for closed orientable 3-manifolds. In
particular, the existing crystallization catalogue $\mathcal
C^{28}$ available in \cite{[L]} is used to obtain upper bounds for
the complexity of closed orientable 3-manifolds triangulated by at
most 28 tetrahedra. The experimental results actually coincide
with the exact values of complexity, for all but three elements.
Moreover, in the case of at most 26 tetrahedra, the exact value of
the complexity is shown to be always directly computable via
crystallization theory.

\end{abstract}

\medskip



\bigskip
\bigskip

\section{\hskip -0.7cm . Introduction}

\medskip

As it is well-known, Matveev's notion of {\it complexity} is based
on the existence, for each compact 3-manifold $M^3$, of a {\it
simple spine}\footnote{According to \cite{[M$_1$]}, a
subpolyhedron $P \subset Int (M^3)$ is said to be a {\it simple
spine} of $M^3$ if the link of each of its points can be embedded
in $\Delta$ (the 1-skeleton of the 3-simplex) and $M^3$ - or $M^3$
minus an open 3-ball, in case $\partial M^3 = \emptyset$ -
collapses to $P$.}: in fact, if  $M^3$ is a compact 3-manifold,
{\it (Matveev's) complexity} \, of  $M^3$ is defined as $$c(M^3)=
min_P c(P),$$ where the minimum is taken over all simple spines
$P$ of $M^3$ and $c(P)$ denotes the number of true
vertices\footnote{Recall that a point of the simple spine $P$ is
said to be a {\it true vertex} if its link is homeomorphic to
$\Delta$.} of the simple spine $P$.

\medskip
In \cite{[C$_4$]} (which is devoted only to the non-orientable
case), a graph-theoretical approach to the computation of
complexity is performed, via another combinatorial theory to
represent 3-manifolds, which makes use of particular edge-coloured
graphs, called {\it crystallizations} (see \cite{[FGG]} or
\cite{[BCG]} for a survey on this representation theory, for
PL-manifolds of arbitrary dimension): the existence of the
crystallization catalogue $\tilde {\mathcal C}^{26}$ (due to
\cite{[C$_2$]}) for closed non-orientable 3-manifolds triangulated
by at most 26 tetrahedra has allowed to complete the existing
classification (due to \cite{[AM]}) of closed non-orientable
3-manifolds up to complexity six.

On the other hand, as already pointed out in \cite{[C$_4$]}, any
crystallization catalogue obviously yields - via the notion of
Gem-Matveev complexity, or GM-complexity, for short, - upper
bounds for the complexity of any involved manifold. Since
complexity and GM-complexity actually turn out to coincide for
each manifold represented by catalogue $\tilde {\mathcal C}^{26}$,
it appears to be an interesting problem to search for classes of
3-manifolds whose complexity can be directly computed via
GM-complexity or, better, to give a characterization of the
classes of 3-manifolds satisfying this property: see
\cite{[C$_4$]} (paragraph 1 - Open Problem).

\medskip

The aim of the present paper is to face the above problem in
the {\it orientable} case, by making use of the existing
crystallization catalogue ${\mathcal C}^{28}$ (due to \cite{[L]})
for closed orientable 3-manifolds triangulated by at most 28
tetrahedra.

For this purpose, a slight modification of the notion of
GM-complexity, involving also non minimal crystallizations, is
taken into account.

Algorithmic computation (easily implemented on computer) directly
yields that, for all but three 3-manifolds involved in ${\mathcal
C}^{28}$, GM-complexity and complexity coincide: see Proposition
\ref{proprietà c_GM in C^(28)}; moreover, if the attention is
restricted to orientable 3-manifolds triangulated by at most 26
tetrahedra, then the exact value of the complexity turns out to be
always directly computable via crystallizations.

\bigskip
\bigskip

\section{\hskip -0.7cm . Crystallizations and GM-complexity}

\medskip

In this section, in order to introduce our graph-theoretical
approach to  the computation of complexity, we briefly recall few
basic concepts of the representation theory of PL-manifolds by
crystallizations. For general PL-topology, Heegaard splittings of
3-manifolds and elementary notions about graphs and embeddings, we
refer to \cite{[HW]}, \cite{[He]} and \cite{[W]} respectively.

Crystallization theory represents PL $n$-manifolds by means of
$(n+1)$\textit{-coloured graphs}, that is, it is a representation
theory which can be used in any dimension. On the other hand,
since this paper concerns only 3-manifolds, the following
definitions and results will be given for $n=3$, although they
mostly hold for each $n\geq 0$.

Moreover, throughout the paper all manifolds will be closed and
connected.

Given a pseudocomplex $K$, triangulating a $3$-manifold
$M$, a \textit{coloration} on $K$ is a labelling of its vertices
by $\Delta_3=\{0,1,2,3\}$, which is injective on each simplex
of $K$.

The dual 1-skeleton of $K$ is a (multi)graph $\Gamma$ embedded in
$|K|=M$; we can define on $\Gamma=(V(\Gamma),E(\Gamma))$ an
\textit{edge-coloration} i.e. a map $\gamma :
E(\Gamma)\to\Delta_3$ in the following way: $\gamma(e)=c\ $ iff
the vertices of the face dual to $e$ are coloured by
$\Delta_3-\{c\}$\footnote{Note that an edge-coloration is
characterized by being injective on each pair of adjacent edges of
the graph.}.

The pair $(\Gamma,\gamma)$ is called a \textit{4-coloured
graph representing M} or simply a \textit{gem=graph encoded
manifold} (see \cite{[L]}).

In the following, to avoid long notations, we will often omit the
edge-coloration, when it is not necessary, and we will simply
write $\Gamma$ instead of $(\Gamma,\gamma)$.

It is easy to see that, starting from $\Gamma$, we can always
reconstruct $K(\Gamma)=K$ and hence the manifold $M$ (see
\cite{[FGG]} and \cite{[BCG]} for more details).

Given $i,j\in\Delta_3$, we denote by $(\Gamma_{i,j},\gamma_{i,j})$
the $3$-coloured graph such that $\Gamma_{i,j}=(V(\Gamma),$
$\gamma^{-1}(\{i,j\}))$ and
$\gamma_{i,j}=\gamma_{|_{\gamma^{-1}(\{i,j\})}}$ i.e. it is
obtained from $\Gamma$ by deleting all edges which are not $i$- or
$j$-coloured; the connected components of $\Gamma_{i,j}$ will be
called \textit{$\{i,j\}$-residues} of $\Gamma$ and their number
will be denoted by $g_{i,j}$.

As a consequence of the definition, a bijection is established
between the set of $\{i,j\}$-residues of $\Gamma$
and the set of 1-simplices of $K(\Gamma)$, whose endpoints
are labelled by $\Delta_3-\{i,j\}$.

Moreover, for each $c\in\Delta_3$, the connected components of the
3-coloured graph $\Gamma_{\hat c}$ obtained from $\Gamma$ by
deleting all $c$-coloured edges, are in bijective correspondence
with the $c$-coloured vertices of $K(\Gamma)$; we will call
$\Gamma$ \textit{contracted} iff $\Gamma_{\hat c}$ is connected
for each $c\in\Delta_3$, i.e. if $K(\Gamma)$ has exactly four
vertices.

A contracted 4-coloured graph representing a $3$-manifold $M$ is called a
\textit{crystallization} of $M$.

Several topological properties of $M$ can be ``read" as
combinatorial properties of any crystallization (or more generally
any gem) $\Gamma$ of $M$: as an example, $M$ is orientable iff
$\Gamma$ is bipartite.

\medskip

Relations among crystallization theory and other classical
representation methods for PL manifolds have been deeply analyzed
(see \cite{[BCG]}; sections 3, 6, 7). In particular, for our
purposes, it is useful to recall the strong connection existing
between crystallizations and Heegaard diagrams.

If $\Gamma$ is a bipartite (resp. non bipartite) crystallization
of a 3-manifold $M,$ for each pair $\alpha, \beta\in \Delta_3,$
let us set $\{\hat\alpha ,\hat\beta\}=\Delta_3-\{\alpha ,\beta\}$
and let $F_{\alpha ,\beta}$ be the orientable (resp. non
orientable) surface of genus $g_{\alpha, \beta}-1=g_{\hat\alpha
,\hat\beta}-1$, obtained from $\Gamma$ by attaching a 2-cell to
each $\{i,j\}$-residue such that $\{i,j\}\neq\{\alpha ,\beta\}$
and $\{i,j\}\neq\{\hat\alpha ,\hat\beta\}$.

It is well-known (see \cite{[FGG]} or \cite{[BCG]}, together with
their references) that a regular embedding\footnote{A cellular
embedding $i$ of a 4-coloured graph $\Gamma$ into a surface is
said to be {\it regular} if there exists a cyclic permutation
$\varepsilon$ of $\Delta_3$ such that the regions of $i$ are
bounded by the images of
$\{\varepsilon_j,\varepsilon_{j+1}\}$-residues of $\Gamma$
($j\in\mathbb Z_4$).} $i_{\alpha, \beta} : \Gamma \to F_{\alpha,
\beta}$  exists. Moreover, if $\mathcal D$ (resp. $\mathcal
D^{\prime}$) is an arbitrarily chosen $\{\alpha, \beta\}$-residue
(resp. $\{\hat\alpha,\hat\beta\}$-residue) of $\Gamma$, the triple
$\mathcal H_{\alpha,\beta,\mathcal D,\mathcal
D^{\prime}}=(F_{\alpha, \beta},{\bf x},{\bf y})$, where {\bf x}
(resp. {\bf y}) is the set of the images of all $\{\alpha,
\beta\}$-residues (resp. $\{\hat \alpha, \hat \beta\}$-residues)
of $\Gamma,$ but $\mathcal D$ (resp. $\mathcal D^{\prime}$), is a
Heegaard diagram of $M$.

Conversely, given a Heegaard diagram $\mathcal H=(F,{\bf x},{\bf
y})$ of $M$ and $\alpha, \beta\in \Delta_3,$ there exists a
construction which, starting from $\mathcal H$ yields a
crystallization $\Gamma$ of $M$ such that $\mathcal H=\mathcal
H_{\alpha,\beta,\mathcal D,\mathcal D^{\prime}}$ for a suitable
choice of $\mathcal D$ and $\mathcal D^{\prime}$ in $\Gamma$ (see
\cite{[G]})\footnote{This correspondence between Heegaard diagrams
and crystallizations allows to prove the coincidence between the
Heegaard genus of $M$ and its {\it regular genus}, a combinatorial
PL-manifold invariant, based on regular embeddings, which is
defined in arbitrary dimension. Interesting results about
classification of PL-manifolds via regular genus may be found, for
example, in \cite{[FG]}, \cite{[CG]}, \cite{[C$_1$]},
\cite{[CM]}}.

Now, let us denote by $\mathcal R_{\mathcal D, \mathcal
D^{\prime}}$ the set of regions of $F_{\alpha, \beta}-({\bf x}\cup
{\bf y})=F_{\alpha, \beta}-i_{\alpha, \beta}((\Gamma_{\alpha,
\beta} - \mathcal D)\cup (\Gamma_{\hat \alpha, \hat \beta} -
\mathcal D^{\prime})).$

\smallskip

\par \noindent
{\bf Definition 3.} Let $M$ be a closed 3-manifold, and let
$\Gamma$ be a crystallization of $M.$  With the above notations,
{\it Gem-Matveev complexity} (or simply {\it GM-complexity}) of
$\Gamma$ is defined as the non-negative integer
$$ c_{GM}(\Gamma) = min \left\{\# V(\Gamma)-\# (V(\mathcal D)\cup V(\mathcal D^{\prime}) \cup V(\Xi)) \ /
\ \alpha, \beta\in \Delta_3, \mathcal D \in \Gamma_{\alpha,\beta},
\mathcal D^{\prime} \in \Gamma_{\hat \alpha, \hat \beta}, \Xi \in
\mathcal R_{\mathcal D, \mathcal D^{\prime}}\right\}$$ while \,
{\it (non-minimal) GM-complexity} \! of $M$ is defined as the
minimum value of GM-complexity, where the minimum is taken over
all crystallizations of $M:$
$$ c^{\prime}_{GM}(M) = min \{c_{GM}(\Gamma) \ / \ \Gamma
\ \text{ is a crystallization \ of \ } M \}$$

\bigskip
As a direct consequence of the definition, (non-minimal)
GM-complexity turns out to be an upper estimation of the manifold
complexity:

\begin{prop} {\rm \bf \cite{[C$_5$]}} \label{relazione
c/cprimeGM} For every closed 3-manifold $M$,
$$c(M) \le c^{\prime}_{GM}(M)$$
\vskip -0,8truecm \ \ \qed \end{prop}

\noindent \textbf{Remark 1.} Definition 3 is a slight modification
(already suggested in \cite{[C$_5$]}) of the previous definition
of {\it $GM$-complexity} of a 3-manifold $M$ denoted by
$c_{GM}(M)$ and originally introduced in \cite{[C$_4$]}. In fact,
$c_{GM}(M)$ is defined as the minimum value of $c_{GM}(\Gamma)$,
too, but the minimum is taken only over \textit{minimal}
crystallizations $\Gamma$ of the manifold, i.e. crystallizations
of $M$ having minimal number of vertices.

In the present paper we will always refer to (non-minimal)
GM-complexity $ c^{\prime}_{GM}(M)$ but, for sake of conciseness,
we will simply write GM-complexity.

\bigskip
\bigskip

It is well known that complexity is additive with respect to the
connected sum of manifolds; GM-complexity can be easily proved to
be subadditive as shown in the following

\begin{prop}\label{subadd} For each pair of closed 3-manifolds
$M_1,M_2,$ the following inequality holds: $$c^{\prime}_{GM}(M_1\#
M_2)\leq c^{\prime}_{GM}(M_1)+c^{\prime}_{GM}(M_2).$$
\end{prop}

\smallskip

\noindent  {\it Hint of the proof.} \  The proof consists
essentially of two steps.

\noindent \underline{Step 1.} \ Let $\Gamma^{(1)},\Gamma^{(2)}$ be
crystallizations of $M_1$ and $M_2$ respectively, such that
$c^{\prime}_{GM}(M_k)=c_{GM}(\Gamma^{(k)})$ ($k=1,2$). With the
notations of Definition 3, for each $k=1,2$, let $\mathcal H_k$ be
the Heegaard diagram associated to $\Gamma^{(k)}$ and $\Xi_k$ the
region of $\mathcal H_k$ realizing $c^\prime_{GM}(\Gamma^{(k)});$
if we denote by $n_k$ the number of vertices of $\mathcal H_k$ and
by $m_k$ the number of vertices of $\Xi_k$, then
$c_{GM}(M_k)=n_k-m_k$. We perform the connected sum of $M_1$ and
$M_2$ with respect to two 3-balls $B_1$ and $B_2$ such that, for
each $k=1,2$, $B_k$ is contained in one of the two handlebodies
defined by $\mathcal H_k$ and intersects the Heegaard surface of
$\mathcal H_k$ in a 2-disc $D_k$ contained in $int\;\Xi_k$. In
this way we obtain a Heegaard diagram $\mathcal H$ of $M_1\# M_2$
having $n_1+n_2$ vertices and containing a region
$\Xi=\Xi_1\#\Xi_2$ with $m_1+m_2$ vertices.

\noindent \underline{Step 2.} \ By applying to the diagram
$\mathcal H$ the construction of \cite{[G]}, a crystallization
$\Gamma$ of $M_1\# M_2$ is obtained, with the property that
$c^{\prime}_{GM}(\Gamma)\leq n_1+n_2-(m_1+m_2)=
c^{\prime}_{GM}(\Gamma^{(1)})+c^{\prime}_{GM}(\Gamma^{(2)})$. \qed

\bigskip

\noindent \textbf{Remark 2.} The additivity of $c^{\prime}_{GM}$
can be proved for the restricted class of manifolds having
GM-complexity coinciding with the complexity. In fact, in this
case, the additivity of the complexity and the above Proposition
yield
$$c^{\prime}_{GM}(M_1\#
M_2)\leq
c^{\prime}_{GM}(M_1)+c^{\prime}_{GM}(M_2)=c(M_1)+c(M_2)=c(M_1\#
M_2)\leq c^{\prime}_{GM}(M_1\# M_2).$$

This result ensures that, as far as we are interested in the
coincidence of GM- and Matveev's complexity, we can restrict our
attention to prime manifolds.

Actually, direct computation proves that the additive property holds for all manifolds
represented by catalogue $\mathcal C^{28}$ (as it is already known for catalogue $\tilde{
\mathcal C}^{26}$).

\bigskip
\bigskip
\bigskip

\section{\hskip -0.7cm . Experimental data from catalogue $\mathcal C^{28}$}

\smallskip

In the literature, a lot of subsequent cataloguing results for
closed orientable irreducible 3-manifolds according to their
complexity exist: in \cite{[M$_1$]} Matveev himself lists all such
manifolds with complexity $c\le 6$; in \cite{[O]} Ovchinnikov
obtains a table for $c=7$ (see also \cite{[M$_2$]} - Appendix 9.3,
where part of Ovchinnikov's table is reproduced); in
\cite{[MP$_1$]} Martelli and Petronio re-obtain via bricks
decomposition the previous results and extend the catalogue up to
complexity 9 (see also
http://www.dm.unipi.it/pages/petronio/public\_html/files/3D/c9/c9\_census.html
for explicit censuses); finally in \cite{[M$_3$]} (see also
\cite{[M$_4$]}) Matveev improves the above classifications by
solving the cases $c=10$ and $c=11$.

Note that closed (orientable and non-orientable\footnote{Burton's
approach allows to include also the case of irreducible and
$\mathbb P^2$-irreducible non-orientable 3-manifolds, which are
further classified up to complexity 7 in \cite{[B$_2$]}.})
3-manifolds up to complexity 6 are also classified by Burton's PhD
thesis, which contains a catalogue of their minimal
triangulations, obtained by face-pairing graphs: see
\cite{[B$_1$]}.

\smallskip

The aim of the present section is to compare, by means of experimental
results, complexity and GM-complexity of orientable 3-manifolds
with ``small" coloured decompositions; for this purpose, the
existing catalogue $\mathcal C^{28}$ of rigid and bipartite
crystallizations with at most 28 vertices (\cite{[L]}) is a basic
tool.

In fact, $\mathcal C^{28}$ yields a catalogue of closed orientable
3-manifolds, ordered by the minimal number of tetrahedra in their
coloured triangulations:

\begin{prop}  {\rm \bf \cite{[L]}} There exist exactly sixty-nine
closed connected prime orientable 3-manifolds, which admit a
coloured triangulation consisting of at most 28 tetrahedra. They
are: the sphere \ $\mathbb S^3$; the orientable $\mathbb
S^2$-bundle over $\mathbb S^1$ (i.e. $\mathbb S^2 \! \times \!
\mathbb S^1$); the six Euclidean orientable 3-manifolds;
twenty-three lens spaces; twenty-one quotients of $\mathbb S^3$ by
the action of their finite (non-cyclic) fundamental groups;
further seventeen topologically undetected orientable
3-manifolds.\footnote{Lins's classification simply identifies
these seventeen 3-manifolds by means of their fundamental groups.
However, five of these groups are given in \cite{[L]} as
semidirect products of $\Bbb Z$ by $\Bbb Z \times \Bbb Z$ induced
by matrices of $GL (2; \Bbb Z)$, and in \cite{[C$_3$]} the
corresponding 3-manifolds are actually proved to be torus bundles
over $\mathbb S^1$.}
\end{prop}

A direct estimation for GM-complexity can be performed for all
manifolds represented in $\mathcal C^{28}$ by means of an easily
implemented computer program\footnote{The C++ program for
GM-complexity computation is available on the Web:
http://cdm.unimo.it/home/matematica/casali.mariarita/DukeIII.htm},
which works as follows:
\begin{itemize}
\item{} given a crystallization $(\Gamma,\gamma)$, let us fix
\begin{itemize}
    \item a partition
    $\{\{\varepsilon_0,\varepsilon_1\},\{\varepsilon_2,\varepsilon_3\}\}$
    of $\Delta_3$;
    \item an $\{\varepsilon_0,\varepsilon_1\}$-residue
    $\mathcal D$ of $\Gamma$;
    \item an $\{\varepsilon_2,\varepsilon_3\}$-residue
    $\mathcal D^\prime$ of $\Gamma$;
    \item a pair of integers $i,j\in\{0,1\}$;
    \item an $\{\varepsilon_i,\varepsilon_{j+2}\}$-residue
    $\Xi_0$ of $\Gamma$;
\end{itemize}

\item{} consider the following subgraph of $\Gamma$,
$$\Xi_1=\left(\bigcup_{\substack{e\in E(\Xi_0)\cap
E(\mathcal D)\\\gamma(e)=\varepsilon_i}}
\Gamma_{\varepsilon_i,\varepsilon_{2+(j+1)mod\,2}}(e)\right)\cup\left(\bigcup_{\substack{f\in
E(\Xi_0)\cap E(\mathcal D^\prime)\\\gamma(f)=\varepsilon_{j+2}}}
\Gamma_{\varepsilon_{(i+1)mod\,2},\varepsilon_{j+2}}(f)\right)$$
where $\Gamma_{a,b}(c)$ denotes the connected component of
$\Gamma_{a,b}$ containing the edge $c$.
 \item{} construct inductively the sequence $\{\Xi_k\}_{k=1,\ldots,m}$, where
 $$\Xi_k=\left(\bigcup_{\substack{e\in E(\Xi_{k-1})\cap
E(\mathcal D)\\
\gamma(e)=\varepsilon_i}}
\Gamma_{\varepsilon_i,\varepsilon_{2+(j+1)mod\,2}}(e)\right)\cup\left(
\bigcup_{\substack{f\in E(\Xi_{k-1})\cap E(\mathcal D^\prime)\\
\gamma(f)=\varepsilon_{j+2}}}
\Gamma_{\varepsilon_{(i+1)mod\,2},\varepsilon_{j+2}}(f)\right)$$
and $m\in\mathbb Z$ is such that $E(\Xi_m)\cap(E(\mathcal D)\cup
E(\mathcal D^\prime))=\emptyset ;$ \item{} compute the number \ \
$\bar c=\# V(\Gamma)-\#(V(\mathcal D)\cup V(\mathcal D^\prime)\cup
V(\Xi_m))$ \ for all possible choices of \ $\Xi_0, i, j, \mathcal
D^\prime, \mathcal D$ \ and for all possible partitions \
$\{\{\varepsilon_0,\varepsilon_1\},\{\varepsilon_2,\varepsilon_3\}\}$
\ of $\Delta_3$; it is very easy to check that the minimal value
assumed by variable $\bar c$
    exactly coincides with $c_{GM}(\Gamma)$;
\item{} for each 3-manifold $M^3$ represented in $\mathcal C^{28}$,
        the above algorithm is applied to the first crystallization $\bar \Gamma$ of
        $M^3$ listed in the catalogue, and then to every
        crystallization $\Gamma \in \mathcal C^{28}$ with
        $K(\Gamma)=M^3;$ the minimal value of their GM-complexities obviously yields an upper estimation for $c^{\prime}_{GM}(M^3)$.
\end{itemize}

\medskip

The obtained results are shown in details in Table 1 of
\cite{[CC-Web]}.  That Table also contains, for each prime
3-manifold involved in $\mathcal C^{28},$  the corresponding
Matveev's description (see \cite{[M$_2$]} - Appendix 9.1 and
Appendix 9.3) and/or the associated Seifert structure: in fact,
for each element of Lins's catalogue, we have also performed the
``translation" into Matveev's notation, which allows a more
efficient topological identification and a direct knowledge of the
complexity. The identifications were usually carried out through
the computation of GM-complexity  and the comparison of homology
groups, and in some cases with the aid of a powerful computer
program for 3-manifold recognition elaborated by Matveev and his
research group and written by V.Tarkaev.\footnote{Computer program
{\it ``Three-manifold Recognizer"} is available on the Web:
http://www.csu.ac.ru/$\sim$trk/}

\medskip

As a consequence, the following improvement of Lins's
classification is obtained, with unambiguous identification of the
encoded 3-manifolds, via JSJ decompositions and fibering
structures: \footnote{In the statement of Proposition \ref{nuova
classificazione C^28}, the following conventions are assumed:
\begin{itemize}
\item[-] for each matrix $A\in GL (2; \Bbb Z)$ with $\det (A)=+1,$ \ $TB(A)= T \times I/A$ \
is the orientable torus bundle over $\mathcal S^1$ with monodromy
induced by $A$;
\item[-] for each matrix $A\in GL (2; \Bbb Z)$ with $\det (A)=-1,$ \ $(K
\tilde \times I) \cup (K \tilde \times I)/A$ \ is the orientable
3-manifold obtained by pasting together, according to $A,$ two
copies of the orientable $I-$bundle over the Klein bottle $K$;
\item[-] $(F,(p_1,q_1), , \dots, (p_k,q_k))$ is the Seifert fibered
manifold with base surface $F$ and $k$ disjoint fibres, having
$(p_i,q_i),$ \ $i=1, \dots, k $ as non-normalized parameters;
\end{itemize}
Moreover, the geometric structures are given according to
\cite{[S]}.}

\begin{prop} \label{nuova
classificazione C^28} The sixty-nine closed connected prime
orientable 3-manifolds which admit a coloured triangulation
consisting of at most 28 tetrahedra are:
\begin{itemize}
\item{} $\mathbb S^3$; \item{}  $\mathbb S^2 \! \times \! \mathbb S^1$; \item{} the six
Euclidean orientable 3-manifolds; \item{} twenty-three lens
spaces; \item{} twenty-one quotients of $\mathbb S^3$ by the
action of their finite (non-cyclic) fundamental groups;
\item{} six (non euclidean) torus bundles $TB(A):$
\begin{itemize}
\item[-] the $Nil$ ones, with complexity 6,  associated to matrices $\begin{pmatrix} -1 & 0 \\
-1 & -1 \end{pmatrix}$ and $\begin{pmatrix}  1 & 0 \\ 1 & 1
\end{pmatrix}$;
\item[-] the $Nil$ ones, with complexity 7,  associated to matrices $\begin{pmatrix} -1 & 0 \\
-2 & -1 \end{pmatrix}$ and $\begin{pmatrix}   1 & -2 \\ 0 & 1
\end{pmatrix}$;
\item[-] the $Sol$ ones, with complexity 7, associated to matrices
$\begin{pmatrix} 0 & 1 \\ -1 & 3 \end{pmatrix}$ and
$\begin{pmatrix}  0 & 1 \\ -1 & -3 \end{pmatrix}$;
\end{itemize}
\item{} two  $Nil$ 3-manifolds of type \ $(K \tilde \times I) \cup (K \tilde \times
I)/A,$ with complexity 6, i.e. the ones associated to matrices
$\begin{pmatrix} 1 & 1 \\ 1 & 0
\end{pmatrix}$ and $\begin{pmatrix} -1 & 0 \\ -1 & 1
\end{pmatrix}$;
\item{} another $Nil$ 3-manifold with complexity 6, i.e. the manifold with Seifert structure $(\mathbb
S^2,(3,2),(3,1),(3,-2))$;
\item{} eight Seifert 3-manifolds with complexity 7:
\begin{itemize}
\item[-] the $Nil$ \, 3-manifold \, $(\mathbb R\mathbb P^2,(2,1),(2,3));$
\item[-] the $SL_2 \mathbb R$ \, 3-manifolds \, $(\mathbb R\mathbb P^2,$
$(2,1),(3,-1)),$ \,  $(\mathbb S^2,$ $ (2,1),$ $(2,1),$ $(2,1),$
$(3,-4)),$
 \, $(\mathbb S^2,(2,1),(3,1),(7,-6)),$
\, $(\mathbb S^2,(2,1),(4,1),(5,-4)),$
 \, $(\mathbb S^2,(3,1),(3,1),(5,-3)),$ \par
 \, $(\mathbb S^2,(3,1),(3,1),(4,-3))$ \, and \,
 \, $(\mathbb S^2,(2,1),(3,1),(7,-5)).$
\end{itemize}
\end{itemize}
\vskip -0,8truecm \ \ \qed
\end{prop}

\bigskip

Experimental data from catalogue $\mathcal C^{28}$ yield
interesting information in order to compare different complexity
notions.

\medskip
First of all, we can consider, together with the complexity, the
so called {\it gem-complexity} of $M^3$, i.e. the non-negative
integer $k(M^3)=\frac {p}{2} -1$, $p$ being the minimum order of a
crystallization of $M^3$: see, for example, \cite{[C$_2$]} -
paragraph 5 or \cite{[C$_4$]} - Remark 1, where the problem of
possible relations between the complexity $c(M^3)$ and
gem-complexity $k(M^3)$ is pointed out.

On one hand, catalogues $\mathcal C^{28}$ and $\tilde{\mathcal
C}^{26}$ allow us to check that, for the first segments of
3-manifold censuses, ``restricted" gem-complexity implies
``restricted" complexity:

\begin{prop}  \label{complessità in C^(28)}
Let $M^3$ be a closed 3-manifold.
\begin{itemize}
\item[(a)] $k(M^3) \le 12 \ \ \ \ \Longrightarrow \ \ \ c(M^3)\le 6;$
\item[(b)] If $M^3$ is assumed to be orientable, then \ \ \ $ k(M^3) \le 13 \ \ \ \ \Longrightarrow \ \ \ c(M^3)\le 7.$
\end{itemize}
\end{prop}

\dimo Statement (b) is a direct consequence of Proposition
\ref{nuova classificazione C^28}, since $c(M^3)\le 7$ holds for
every manifold $M^3$ encoded by elements of $\mathcal C^{28}.$

Statement (a) concerns both orientable and non-orientable
3-manifolds. In the orientable case, it follows from
identification results contained in the first part of Table 1 of
\cite{[CC-Web]}: in fact, $k(M^3) \le 12$ implies the existence of
a rigid cristallization $\bar \Gamma \in \mathcal C^{28}$, with
$\#V(\bar \Gamma)\le 26$, representing $M^3$, and this immediately
yields $c(M^3)\le 6$ (as it may be seen in the last column of
Table 1 itself). On the other hand, in the non-orientable case,
statement (a) is a direct consequence of results contained in
\cite{[C$_2$]} and \cite{[C$_4$]} (see also \cite{[C$_4$]} -
Remark 1, where the set of irreducible and $\mathcal
P^2$-irreducible non-orientable 3-manifolds up to complexity $c=6$
is proved to coincide exactly with the set of such manifolds up to
gem-complexity $k=12$).
 \qed

\bigskip

On the other hand, for all manifolds in the catalogue $\mathcal
C^{28}$ ``restricted" complexity implies ``restricted"
gem-complexity, too. More precisely, we can state:

\begin{prop}  \label{gem_complexity-complexity in C^(28)}
Let $M^3$ be a closed orientable 3-manifold with complexity
$c(M^3)=c.$  \  If \, $0 \le c\le 4,$ \, then  \ $k(M^3) \le
5+2c.$
\end{prop}

\dimo It is well-known, within crystallization theory, that
$\mathbb S^1 \times \mathbb S^2$ admits a (non-rigid) order eight
crystallization; hence, $k(\mathbb S^1 \times \mathbb S^2)= 3$.
This fact, together with a direct comparison between Table 1 of
\cite{[CC-Web]} and the tables of \cite{[M$_2$]} - Appendix 9.1
allows to  state that all closed orientable 3-manifolds with
complexity 0 (resp. 1) (resp. 2) (resp. 3) (resp. 4) admit a gem
with at most 12 (resp. 16) (resp. 20) (resp. 24) (resp. 28)
vertices. Hence, the corresponding gem-complexities satisfy the
claimed inequality. \qed

The above results naturally suggest the following

\medskip \noindent {\bf Conjecture:}
\ $ k(M^3) \le 5 +2 c(M^3)$ \ \ for any closed orientable
3-manifold $M^3.$

\bigskip
\bigskip

Moreover, experimental data concerning GM-complexity estimation
for closed orientable 3-manifolds represented by the
crystallization catalogue $\mathcal C^{28}$ - appearing in the
fifth column of Table 1 of \cite{[CC-Web]}, - allow us to prove directly the following properties, and therefore to establish a
comparison between GM-complexity and complexity. Note that, for
sake of notational simplicity, $M^3 \in \mathcal C^{2p}$ ($p\in
\mathbb Z$) is written in order to indicate
a manifold $M^3$ which admits a rigid crystallization belonging to
the catalogue $\mathcal C^{2p}$ of all rigid bipartite
crystallizations with order $\le 2p.$

\begin{prop}  \label{proprietà c_GM in C^(28)}
\par \noindent
\begin{itemize}
\item[(a)] $c^{\prime}_{GM}(M^3) \le c(M^3) +1 $ \ \ \ \ $ \forall M^3 \in \mathcal
C^{28};$
\item[(b)]
$ c^{\prime}_{GM}(M^3) = c(M^3)$ \ \ \ \ \ \ \ \ \  $\forall M^3
\in \mathcal C^{26};$
\item[(c)] $c^{\prime}_{GM}(M^3) = c(M^3)$ \ \ \ \ \ \ \ \ \ \ \ \ \ $
\forall M^3 \in \mathcal C^{28}- \Big\{ (\mathbb R\mathbb P^2,$
$(2,1),(3,-1)),$ \ \ $(\mathbb S^2,(3,1),(3,1),(5,-3)),$ \  \par \
\hskip 7.2truecm $(\mathbb S^2, (2,1), (2,1),(2,1),(3,-4))
\Big\}.$
\end{itemize}
\vskip -0,6truecm \ \ \qed
\end{prop}

\bigskip
\noindent  {\bf Remark 3.} It is an open problem to compute the
values (belonging to the set $\{7,8\}$) of
$c^{\prime}_{GM}\left((\mathbb R\mathbb P^2,(2,1),(3,-1))\right),$
\ $c^{\prime}_{GM}\left((\mathbb S^2,(3,1),(3,1),(5,-3))\right),$
\ $c^{\prime}_{GM}\left((\mathbb
S^2,(2,1),(2,1),(2,1),(3,-4))\right).$

\bigskip
\noindent  {\bf Remark 4.} For all manifolds, but one, encoded in
catalogue $\mathcal C^{28}$ and whose GM-complexity and complexity
coincide, GM-complexity is realized by a minimal cristallization
in the sense of gems (according to the original definition of
GM-complexity $c_{GM}(M)$, introduced in \cite{[C$_4$]}). More
precisely, if $\mathbb S^3/G$ denotes the quotient space of
$\mathbb S^3$ by the action of group $G$ and $P_{24}=<x,y\, |\,
x^2=(xy)^3=y^3, x^4=1>$, the following result holds:
\par \noindent $ \forall M^3 \in \mathcal C^{28}$ so that
$c^{\prime}_{GM}(M^3) = c(M^3)$, $M^3 \ne \mathbb S^3/(P_{24}
\times \mathbb Z_5),$ \ then \ $c^{\prime}_{GM}(M^3)=
c_{GM}(\Gamma)$, with $\#V(\Gamma)\le \#V(\Gamma^{\prime})$,
$\forall \Gamma^{\prime}$ representing $M^3$.


\begin{thebibliography}{[BCG]}


\bibitem{[AM]} G.Amendola - B.Martelli {\em Non-orientable 3-manifolds
of small complexity}, Topology Appl. {\bf 133} (2003), 157-178.

\bibitem{[B$_1$]} B.A. Burton, {\em Minimal triangulations and normal surfaces}, PhD Thesis, University of Melbourne (Australia), May 2003,
available from the Web page
http://regina.sourceforge.net/data.html

\bibitem{[B$_2$]} B.A. Burton, {\em Structures of small
closed non-orientable 3-manifold triangulations}, Math.GT/0311113.

\bibitem{[BCG]} P.Bandieri - M.R.Casali - C.Gagliardi,
{\em Representing manifolds by crystallization theory:
foundations, improvements and related results}, Atti Sem. Mat.
Fis. Univ. Modena {\bf Suppl. 49} (2001), 283-337.

\bibitem{[C$_1$]} M.R.Casali, {\em Classifying PL 5-manifolds by regular genus:
the boundary  case}, Canadian J. Math. {\bf 49}  (1997), 193-211.

\bibitem{[C$_2$]} M.R.Casali, {\em
Classification of non-orientable 3-manifolds admitting
decompositions into $\le$ 26 coloured tetrahedra}, Acta Appl.
Math. {\bf 54} (1999), 75-97.

\bibitem{[C$_3$]} M.R.Casali, {\em
Representing and recognizing torus bundles over $\mathbb S^1$},
Boletin de la Sociedad Matematica Mexicana (special issue in honor
of Fico), {\bf 10 (3)}  (2004), to appear.

\bibitem{[C$_4$]} M.R.Casali, {\em Computing Matveev's complexity of non-orientable
3-manifolds via crystallization theory}, Topology and its
Applications {\bf 144} (2004), 201-209.

\bibitem{[C$_5$]} M.R.Casali, {\em Estimating Matveev's complexity via crystallization theory},
to appear.

\bibitem{[CC-Web]} M.R.Casali - P. Cristofori, {\em Archives of closed 3-manifolds with low gem-complexity}, available from the Web page
http://cdm.unimo.it/home/matematica/casali.mariarita/DukeIII.htm

\bibitem{[CG]} M.R.Casali - C.Gagliardi, {\em Classifying
PL 5-manifolds up to regular genus seven}, Proc. Amer. Math. Soc.
{\bf 120 (1)} (1994), 275-283.

\bibitem{[CM]} M.R.Casali - L.Malagoli, {\em Handle-decompositions of PL
4-manifolds}, Cahiers de Topologie et Geom. Diff. Cat. {\bf 38}
(1997), 141-160.

\bibitem{[FG]} M.Ferri - C.Gagliardi, {\em The only
genus zero n-manifold is $\mathbb S^{\text {\it n}}\, $},  Proc.
Amer. Math. Soc. {\bf 85} (1982), 638-642.

\bibitem{[FGG]} M.Ferri - C.Gagliardi - L.Grasselli, {\em A
graph-theoretical representation of PL-manifolds. A survey on
crystallizations},   Aequationes Math. {\bf 31} (1986), 121-141.

\bibitem{[G]} C.Gagliardi, {\em Extending the concept of genus
to dimension $n$}, Proc. Amer. Math. Soc. {\bf 81} (1981),
473-481.

\bibitem{[He]} J.Hempel, {\em 3-manifolds},  Annals of Math. Studies, {\bf 86},
Princeton Univ. Press, 1976.

\bibitem{[HW]} P.J.Hilton - S.Wylie, {\em An introduction to
algebraic topology - Homology theory},  Cambridge Univ. Press,
1960.

\bibitem{[L]} S.Lins, {\em Gems, computers and attractors for
3-manifolds}, Knots and Everything {\bf 5}, World Scientific,
1995.

\bibitem{[MP$_1$]} B.Martelli - C.Petronio, {\em Three-manifolds
having complexity at most 9}, Experimental Mathematics {\bf 10
(2)} (2001), 207-236.

\bibitem{[MP$_2$]} B.Martelli - C.Petronio, {\em Census 7},
Table of closed orientable irreducible 3-manifolds having
complexity 7, available from the Web page
http://www.dm.unipi.it/pages/petronio/public\_html/files/3D/c9/c9\_census.html

\bibitem{[M$_1$]} S.Matveev, {\em Complexity theory of
three-dimensional manifolds}, Acta Applicandae Math. {\bf 19}
(1990), 101-130.

\bibitem{[M$_2$]} S.Matveev, {\em  Algorithmic topology and  classification of 3-manifolds},
Algorithms and Computation in Mathematics {\bf 9}, Springer, 2003.

\bibitem{[M$_3$]} S.Matveev, {\em Recognition and tabulation of three-dimensional manifolds},
Doklady RAS {\bf 400(1)}(2005), 26-28 (Russian;  English trans. in
Doklady Mathematics, {\bf 71} (2005), 20-22).

\bibitem{[M$_4$]} S.Matveev, {\em Tabulation of 3-manifolds},
Uspekhi Mt. Nauk. {\bf 60(4)}(2005), 97-122 (Russian; English
trans. in Russian Math. Surveys {\bf 60(4)}(2005), 673-698).


\bibitem{[O]} M.A. Ovckinnikov, {\em The table of 3-manifolds of complexity
7}, Preprint Chelyabinsk State University, 1997.

\bibitem{[S]} P.Scott, {\em The geometries of 3-manifolds},  Bull. London Math. Soc.  {\bf 15} (1983),
401-487.

\bibitem{[W]} A.T.White, {\em Graphs, groups and surfaces}, North
Holland, 1973.


\end{thebibliography}
\end{document}